%Paper: alg-geom/9312001
%From: David Cox <dac@cs.amherst.edu>
%Date: Thu, 2 Dec 1993 21:33:58 -0500
%Date (revised): Wed, 8 Dec 1993 16:04:36 -0500

% This is the tex file for ``The Functor of a Smooth Toric Variety''
%
% This file is written in plain TeX, but it assumes that you have
% access to the AMS-TeX fonts.
%
% If you have any problems with the blackboard bold fonts, delete all
% of the lines in the section titled ``Blackboard Roman Letters'' and
% replace it with the single line:
%
% \def{\bb#1{{\bf\relax#1}}
%
% If you have any problems texing the file, contact David Cox.  My
% email address is dac@cs.amherst.edu
%
% Here are the macros used in the text.  We begin with the formatting
% macros:

\magnification=\magstep1
\baselineskip=14pt
\abovedisplayskip=10pt plus 2pt minus 2pt
\belowdisplayskip=10pt plus 2pt minus 2pt
\font\bigbf=cmbx10 scaled \magstep1 % large font

% Blackboard Roman letters

\font\tenmsy=msym10
\font\sevenmsy=msym7
\font\fivemsy=msym5
\newfam\msyfam
\textfont\msyfam=\tenmsy
\scriptfont\msyfam=\sevenmsy
\scriptscriptfont\msyfam=\fivemsy
\def\bb#1{{\fam\msyfam\relax#1}}

% Commonly used letters (Greek and Blackboard Bold)

\def\a{\alpha}
\def\A{{\bb A}}
\def\P{{\bb P}}
\def\r{{\bb R}}
\def\z{{\bb Z}}

% Other useful macros

\def\d{$$}
\def\prf{\noindent {\bf Proof.\ \ }}
\def\square{\smash{\raise8pt\hbox{\leavevmode\hbox{\vrule\vtop{\vbox
	{\hrule width9pt}\kern9pt\hrule width9pt}\vrule}}}}
\def\qed{\hfill\square}
\def\ref{\hangindent=\parindent \hangafter=1 \noindent}
\def\mapname#1{\ \smash{\mathop{\longrightarrow}\limits^{#1}}\ }

% Macros special to this paper:

\def\pic{{\rm Pic}(X)}
\def\oy{{\cal O}_Y}
\def\ox{{\cal O}_X}
\def\ou{{\cal O}_U}
\def\oU{{\cal O}_{{\cal U}}}
\def\oP{{\cal O}_{\P^m_k}}
\def\xcoll{{$\Delta$-collection}}
\def\hom{{\rm Hom}_k}

% Now comes the text of the paper:

\centerline{\bigbf The Functor of a Smooth Toric Variety}
\bigskip

\centerline{David A. Cox}
\centerline{Department of Mathematics and Computer Science}
\centerline{Amherst College}
\centerline{Amherst, MA 01002}
\centerline{dac@cs.amherst.edu}
\bigskip

	In EGA I [3], projective space $\P_k^n$ is described as
the variety representing the functor
\d Y \mapsto \{\hbox{line bundle quotients of $\oy^{n+1}$}\}
\ .\leqno(1)
\d
This is easy to prove since a surjection $\oy^{n+1} \to L$ gives $n+1$
sections of $L$ which don't vanish simultaneously and hence determine
a map $Y \to \P^n_k$.  The goal of this paper is to generalize this
representation to the case of an arbitrary smooth toric variety.

	We will work with schemes over an algebraically closed field
$k$ of characteristic zero, and we will fix a smooth $n$-dimensional
toric variety $X$ determined by the fan $\Delta$ in $N_\r = \r^n$.  As
usual, $M$ denotes the dual lattice of $N$ and $\Delta(1)$ denotes the
1-dimensional cones of $\Delta$.  We will use $\sum_\rho$ to mean
$\sum_{\rho \in \Delta(1)}$, and similarly for $\otimes_\rho$.  Each
$\rho \in \Delta(1)$ determines a divisor $D_\rho \subset X$ and a
generator $n_\rho$ of $\rho\cap N$.  Finally, let $\Delta_{max}$
denote the maximal cones in $\Delta$ (i.e., those which are not proper
faces of cones in $\Delta$).
\bigskip

\noindent {\bf \S1. $\Delta$-Collections and Functors.} If a fan
$\Delta$ determines a smooth toric variety $X$, then we can generalize
the data in (1) as follows:

\proclaim Definition 1.1. Given a scheme $Y$ over $k$, a {\bf
$\Delta$-collection on $Y$} consists of line bundles $L_\rho$ and
sections $u_\rho \in H^0(Y,L_\rho)$, indexed by $\rho \in \Delta(1)$,
and isomorphisms $c_m : \otimes_\rho L_\rho^{\otimes\langle
m,n_\rho\rangle} \simeq \oy$, indexed by $m \in M$, such that:
\vskip0pt
\item{(i)} (Compatibility) $c_m\otimes c_{m'} = c_{m+m'}$ for all
$m,m' \in M$.
\item{(ii)} (Nondegeneracy) $u_\rho \in H^0(Y,L_\rho)$ gives
$u_\rho : \oy \to L_\rho$, which induces $u_\rho^* :L_\rho^{-1}
\to \oy$.  Then the map $\sum_{\sigma \in
\Delta_{max}} \otimes_{\rho \not\subset \sigma} u_\rho^* :
\bigoplus_{\sigma \in \Delta_{max}} \otimes_{\rho \not\subset \sigma}
L_\rho^{-1} \to \oy$ is a surjection.
\vskip0pt

	A $\Delta$-collection on $Y$ is written
$(L_\rho,u_\rho,c_m)$.  The compatibility condition on the
isomorphsims $c_m$ implies that $\sum_\rho [L_\rho]\otimes n_\rho = 0$
in ${\rm Pic}(Y)\otimes_\z N$.  However, the triviality of this sum is
not sufficient: data of the \xcoll\ includes an explicit choice of
trivialization (the $c_m$'s), which is not unique.  The examples given
below will show why this is needed.  Finally, when dealing with a
toric variety $X$ determined by $\Delta$, we will sometimes speak of
$X$-collections rather than \xcoll s.

	We get a canonical \xcoll\ on the toric variety $X$ as
follows.  For each $\rho$, the divisor $D_\rho$ gives a line bundle
$\ox(D_\rho)$ since $X$ is smooth.  Furthermore, the natural inclusion
$\ox \subset \ox(D_\rho)$ corresponds to a global section $\iota_\rho
\in H^0(X,\ox(D_\rho))$.  Finally, given $m \in M$, the character
$\chi^m$ is a rational function on $X$ such that ${\rm div}(\chi^m) =
\sum_\rho \langle m,n_\rho\rangle D_\rho$.  Thus we get an isomorphism
of sheaves
\d c_{\chi^m} : \otimes_\rho \ox(D_\rho)^{\otimes \langle m,n_\rho\rangle}
\simeq \ox\big(\textstyle{\sum_\rho} \langle m,n_\rho\rangle
D_\rho\big) \simeq \ox\ ,
\d
where the second isomorphism is induced by $\chi^m$.

\proclaim Lemma 1.1. $(\ox(D_\rho),\iota_\rho,c_{\chi^m})$ is a
\xcoll\ on $X$.

\prf Compatibility is trivial since $\chi^{m+m'} = \chi^m\chi^{m'}$
for $m,m' \in M$.  To prove nondegeneracy, take $x \in X$.  Since $X =
\cup_{\sigma \in \Delta_{max}} X_\sigma$, where $X_\sigma$ is the
affine toric variety determined by $\sigma$, we have $x \in X_\sigma$
for some $\sigma \in \Delta_{max}$.  Since $X - X_\sigma = \sum_{\rho
\not\subset \sigma} D_\rho$, we see that $\otimes_{\rho \not\subset
\sigma} \iota_\rho^* : \otimes_{\rho \not\subset \sigma}
\ox(D_\rho)^{-1} = \ox\big(-\sum_{\rho \not\subset \sigma} D_\rho\big)
\to \ox$ is an isomorphism at $x$, and nondegeneracy follows. \qed
\medskip

	The \xcoll\ $(\ox(D_\rho),\iota_\rho,c_{\chi^m})$ will be called
the {\it universal \xcoll\/}.  This terminology will be justified below.

\proclaim Definition 1.2. An {\bf equivalence} $(L_\rho,u_\rho,c_m) \sim
(L'_\rho,u'_\rho,c'_m)$ of $\Delta$-collections on $Y$ consists of
isomorphisms $\gamma_\rho : L_\rho \simeq L'_\rho$ which carry
$u_\rho$ to $u'_\rho$ and $c_m$ to $c'_m$.

	To better understand these definitions, let's look at some
examples:
\medskip

\noindent {\bf Example 1.1.} Let $X = \A^n_k$, where the $n_\rho$'s are
the standard basis $e_1,\dots,e_n$ of $N$.  Let $e^1,\dots,e^n$ be the
dual basis of $M$.  Now suppose we have a $\A^n_k$-collection
$(L_i,u_i,c_m)$ on $Y$ (we write $L_i$ instead of $L_{e_i}$ for $1\le
i \le n$).  Then $c_{e^i}$ is an isomorphism $c_{e^i} : L_i \simeq
\oy$.  This maps $u_i$ to $v_i \in H^0(Y,\oy)$, and one can check that
setting $\gamma_i = c_{e^i}$ in Definition 1.2 gives an equivalence
$(L_i,u_i,c_m) \sim (\oy,v_i,1)$.  Furthermore, $(\oy,v_i,1) \sim
(\oy,v'_i,1)$ if and only if $v_i = v'_i$ for all $i$.  Since the
nondegeneracy condition is vacuous in this case, we see that
equivalence classes of $\A^n_k$-collections on $Y$ correspond exactly
to $n$-tuples in $H^0(Y,\oy)$.  As is well-known, such $n$-tuples are
classified by morphisms $Y \to \A^n_k$.
\medskip

\noindent{\bf Example 1.2.} Let $X = \P^n_k$, where the $n_\rho$'s are
$e_1,\dots,e_n$ and $e_0 = -\sum_{i=0}^n e_i$.  Now let
$(L_i,u_i,c_m)$ be a $\P^n_k$-collection on $Y$ (where $0 \le i \le
n$).  Here, $c_{e^i}$ is an isomorphism $c_{e^i} : L_i\otimes L_0^{-1}
\simeq \oy$.  This induces $\gamma_i : L_i \simeq L_0$ which takes
$u_i$ to $v_i \in H^0(Y,L_0)$.  One can check that the $\gamma_i$'s
give an equivalence $(L_i,u_i,c_m) \sim (L_0,v_i,1)$.  Furthermore,
the nondegeneracy condition means that the map $\sum_i v_i^*:
\bigoplus_{i=0}^n L_0^{-1} \to \oy$ is surjective (the maximal cones
are indexed by $0,...,n$).  Tensoring with $L_0$, this is equivalent
to the surjectivity of $\sum_i v_i : \oy^{n+1} \to L_0$.  It follows
that equivalence classes of $\P^n_k$-collections on $Y$ correspond
exactly to line bundle quotients of $\oy^{n+1}$, and hence, by (1),
are classified by morphisms $Y \to \P^n_k$.
\medskip

\noindent{\bf Example 1.3.} Let $X = {\bb G}_m^n$.  In this case,
there are no $n_\rho$'s, and for $m \in M$, $\otimes_\rho
L_\rho^{\otimes\langle m,n_\rho\rangle}$ reduces to $\oy$, so that a
${\bb G}_m^n$-collection on $Y$ consists of $c_m : \oy \simeq \oy$.
The notions of equivalence and nondegeneracy are vacuous in this case,
so that equivalence classes of ${\bb G}_m^n$-collections on $Y$
correspond to homomorphisms $M \to H^0(Y,\oy^*)$.  Such homomorphisms
are classified by morphisms $Y \to {\rm Hom}_\z(M,{\bb G}_m) = {\bb
G}_m^n$.
\medskip

	Returning to the general case, it is easy to see that the
pull-back of a \xcoll\ is again a \xcoll.  Thus we get a functor $C_\Delta :
k\hbox{-Schemes}^\circ \to {\rm Sets}$ defined by
\d C_\Delta(Y) = \{\hbox{all \xcoll s $(L_\rho,u_\rho,c_m)$ on $Y$}\}/\sim\ .
\d
Furthermore, the universal \xcoll\ $(\ox(D_\rho),\iota_\rho,c_{\chi^m})$
gives a natural transformation
\d \hom(Y,X) \to C_\Delta(Y)
\d
by sending $f:Y \to X$ to the pull-back of
$(\ox(D_\rho),\iota_\rho,c_{\chi^m})$ by $f$.  The main result of this
paper is the following theorem:

\proclaim Theorem 1.1.  If $X$ is a smooth toric variety, then the
above map $\hom(Y,X) \to C_\Delta(Y)$ is a bijection for all
$k$-schemes $Y$.  Thus the toric variety $X$ represents the functor
$C_\Delta$.

\prf First assume that the $n_\rho$'s span $N_\r$.  In this case, we
know by [2] that $X$ is a geometric quotient $(\A_k^{\Delta(1)} -
Z)/G$, where $\A_k^{\Delta(1)} = {\rm Spec}(k[x_\rho])$, $Z$ is
defined by the vanishing of $\prod_{\rho \not\subset \sigma} x_\rho$
for $\sigma \in \Delta_{max}$, and $G = {\rm Hom}_\z(\pic,{\bb G}_m)$.

	We will construct an inverse map $C_\Delta(Y) \to \hom(Y,X)$.
Let $(L_\rho,u_\rho,c_m)$ be a \xcoll\ on $Y$, and let $U \subset Y$
be an open subset such that the $L_\rho$ are trivial on $U$.  If we
choose isomorphisms $\gamma_\rho : L_{\rho|U} \to \ou$, then we get an
equivalence $(L_{\rho|U},u_{\rho|U},c_{m|U}) \sim (\ou,v_\rho,c'_m)$,
where $v_\rho \in H^0(U,\ou)$ and $c'_m: \ou \simeq \ou$ can be
regarded as a homomorphism $c' : M \to H^0(U,\ou^*)$.

	Since the $n_\rho$'s span $N_\r$ and $X$ is smooth, we have an
exact sequence
\d 0 \longrightarrow M \mapname{\a} \z^{\Delta(1)} \longrightarrow
\pic \longrightarrow 0 \leqno(2)
\d
where $\a$ is defined by $m \mapsto (\langle m,n_\rho\rangle)$.  Since
$\pic$ is torsion free, the above map $c' : M \to H^0(U,\ou^*)$
extends to $\tilde c':\z^{\Delta(1)} \to H^0(U,\ou^*)$, which means
that there are $\lambda_\rho \in H^0(U,\ou^*)$ such that $c'_m =
\prod_\rho \lambda_\rho^{\langle m,n_\rho\rangle}$ for all $m \in M$.
Then the isomorphisms $\lambda_\rho : \ou \simeq \ou$ give an
equivalence $(\ou,v_\rho,c'_m) \sim (\ou,w_\rho,1)$, where $w_\rho =
\lambda_\rho v_\rho$.

	Now define $\tilde f_U : U \to \A^{\Delta(1)}_k$ by $\tilde
f_U(x) = (w_\rho(x))$.  The nondegeneracy condition implies that
$\tilde f_U(x) \notin Z$, so that composing with the quotient map $\pi
: \A^{\Delta(1)} - Z \to X$ gives $f_U =\pi\circ\tilde f_U: U \to X$.

	We need to see how the choices made in the above construction
affect the map $f_U$.  A different set of choices would lead to a
\xcoll\ $(\ou,w'_\rho,1) \sim(\ou,w_\rho,1)$.  This equivalence is
given by $\lambda_\rho \in H^0(U,\ou)$ such that $w'_\rho =
\lambda_\rho w_\rho$ and $\prod_\rho \lambda_\rho^{\langle
m,n_\rho\rangle} = 1$ for all $m \in M$ (because the $\lambda_\rho$'s
must preserve the trivializations $1 : \otimes \ou^{\langle
m,n_\rho\rangle} \simeq \ou$).  It follows from (2) that we get a
homomorphism $g: \pic \to H^0(U,\ou^*)$ such that $g([D_\rho]) =
\lambda_\rho$ for all $\rho$.  If we evaluate this at a closed point
$x \in U$, we get an element $g_x \in G = {\rm Hom}_\z(\pic,{\bb
G}_m)$.  Then the points $(w_\rho(x))$ and $(w'_\rho(x)) =
(\lambda_\rho(x) w_\rho(x))$ are related by $g_x$ and hence give the
same point in $X$.

	This shows that $f_U : U \to X$ depends only on the
equivalence class of $(L_\rho,u_\rho,c_m)$.  From here, it follows
easily that the $f_U$ patch together to give a morphism $f : Y \to X$.

	It remains to show that this map is the inverse of the map
$\hom(Y,X) \to C_\Delta(Y)$ obtained by pulling back the univeral \xcoll\
$(\ox(D_\rho), \iota_\rho, c_{\chi^m})$.  First suppose that
$(L_\rho,u_\rho,c_m)$ on $Y$ determines $f:Y\to X$.  We need to show
that
\d (L_\rho,u_\rho,c_m) \sim f^*(\ox(D_\rho),\iota_\rho,c_{\chi^m})\ .
\leqno(3)
\d
An easy argument shows that the natural map $C_\Delta(Y) \to \prod_{\a
\in A} C_\Delta(U_\a)$ is injective whenever $\{U_\a\}_{\a \in A}$ is
an open cover of $Y$.  Thus it suffices to prove (3) on an open set $U
\subset Y$ where each $L_\rho$ is trivial on $U$.  On such a $U$, we
know that $(L_{\rho|U},u_{\rho|U}, c_{m|U}) \sim (\ou,w_\rho,1)$ and
$f = \pi \circ \tilde f$, where $\tilde f(x) = (w_\rho(x))$.  We first
observe that
\d \pi^*(\ox(D_\rho),\iota_\rho,c_{\chi^m}) \sim (\oU,x_\rho,1)\ ,
\leqno(4)
\d
where $\pi : {\cal U} \to X$ is as above.  To prove (4), note that
$\pi^*(\ox(D_\rho)) = \oU({\rm div}(x_\rho))$, so that multiplication
by $x_\rho$ gives an isomorphism $\oU({\rm div}(x_\rho)) \simeq \oU$.
Hence $\oU \subset \oU({\rm div}(x_\rho)) \simeq \oU$ is
multiplication by $x_\rho$, and since $\chi^m\circ \pi = \prod_\rho
x_\rho^{\langle m,n_\rho\rangle}$, (4) follows immediately.  Then,
returning to $f = \pi\circ\tilde f$, we conclude from (4) that
\d f^*(\ox(D_\rho),\iota_\rho,c_{\chi^m}) \sim \tilde f^*(\oU,x_\rho,1)
= (\ou,w_\rho,1)\leqno(5)
\d
since $\tilde f(x) = (w_\rho(x))$, and (3) follows.

	Finally, suppose we have $f:Y \to X$.  This gives
$f^*(\ox(D_\rho), \iota_\rho,c_{\chi^m})$, which in turn determines
$f':Y\to X$.  We need to show that $f' = f$.  First suppose that $f$
factors $f = \pi\circ\tilde f$ for some map $\tilde f : Y \to {\cal
U}$.  Then $\tilde f$ can be written $\tilde f(x) = (w_\rho(x))$,
where $w_\rho \in H^0(Y,\oy)$.  From (5) and the construction of $f'$,
it follows immediately that $f' = f$.  In the general case, note that
$G$ acts freely on ${\cal U}$ since $X$ is smooth (this is easy to
prove), so that $\pi : {\cal U} \to X$ is smooth.  Then standard
results about smoothness imply that $f:Y\to X$ lifts locally to ${\cal
U}$ in the \'etale topology.  Since $\hom(-,X)$ is a sheaf in the
\'etale topology on $Y$, we obtain $f'=f$, and the theorem is proved
in the case when the $n_\rho$'s span $N_\r$.

	We next study what happens when the $n_\rho$'s don't span
$N_\r$.  Let $N_1 = N\cap{\rm Span}_\r(n_\rho)$.  The fan $\Delta$ can
be regarded as a fan $\Delta_1$ in $N_1$, which gives a smooth toric
variety $X_1$ of dimension $d = {\rm rank}(N_1)$.  The inclusion $N_1
\subset N$ induces an inclusion $X_1 \subset X$, and the projection $N
\to N/N_1$ induces a surjection $X \to T_1 = {\rm Hom}_\z(N_1^\perp,
{\bb G}_m) \simeq {\bb G}_m^{n-d}$, where $N_1^\perp = {\rm
Hom}_\z(N/N_1,\z) \subset M$ is the annihilator of $N_1$.

	Since $N/N_1$ is torsion free, we can write $N = N_1\oplus
N_2$ for some complement $N_2 \subset N$.  Then $\Delta$ is the
product fan $\Delta_1\times\{0\}$, which implies that $X$ is
(noncanonically) the product $X_1\times_k T_1$.  If $M_1$ is the dual
of $N_1$, then the projection $N \to N_1$ determines an inclusion $\a
: M_1 \to M$ such that $M = \a(M_1)\oplus N_1^\perp$.

	Now suppose that $(L_\rho,u_\rho,c_m)$ is a \xcoll\ on $Y$.
Then, for every $m \in N_1^\perp$, we have $\langle m, n_\rho\rangle =
0$ for all $\rho$.  Thus $c_m$ is an isomorphism $c_m : \oy \simeq
\oy$, which gives a homomorphism $N_1^\perp \to H^0(Y,\oy^*)$.
Since this map depends only on the equivalence class of $(L_\rho,
u_\rho,c_m)$, we have a natural transformation
\d C_\Delta(Y) \longrightarrow {\rm Hom}_\z(N_1^\perp,H^0(Y,\oy^*))\ .
\d
Further, if we define $c_{m_1}^\a = c_{\a(m_1)}$ for $m_1 \in M_1$,
then $(L_\rho, u_\rho, c^\a_{m_1})$ is a $\Delta_1$-collection on $Y$,
and it follows easily that we have a natural transformation
\d C_\Delta(Y) \longrightarrow C_{\Delta_1}(Y)\ .
\d
Combining these maps, we obtain
\d C_\Delta(Y) \longrightarrow C_{\Delta_1}(Y)\times{\rm
Hom}_\z(N_1^\perp,H^0(Y,\oy^*))\ . \leqno(6)
\d
Since $M = \a(M)\oplus N_1^\perp$, it is straightforward to show that
the map (6) is a bijection.

	Now consider the following diagram:
\d \matrix{\hom(Y,X) & \longrightarrow & C_\Delta(Y)\cr
\downarrow && \downarrow \cr
\hom(Y,X_1)\times\hom(Y,T_1) & \longrightarrow &
C_{\Delta_1}(Y)\times{\rm Hom}_\z(N_1^\perp,H^0(Y,\oy^*))\cr}
\d
The vertical maps come from (6) and $X \simeq X_1\times T_1$, and note
that both are bijections.  The map on the bottom is the product of the
bijections $\hom(Y,X_1) \simeq C_{\Delta_1}(Y)$ (since the $n_\rho$'s
span $(N_1)_\r$) and $\hom(Y,T_1) \simeq {\rm
Hom}_\z(N_1^\perp,H^0(Y,\oy^*))$ (since $T_1 = {\rm Hom}_\z(N_1^\perp,
{\bb G}_m)$).

	It follows that the map on top will be a bijection (and the
theorem will be proved) provided the diagram commutes.  By general
nonsense, we only have to prove commutivity for $1_X \in \hom(X,X)$.
Going down and over, $1_X$ maps to $(\pi_1^*({\cal
O}_{X_1}(D_\rho),\iota_\rho,c_{\chi^{m_1}}),\phi)$, where $\pi_1 : X
\to X_1$ is the projection and $\phi:N_1^\perp \to H^0(X,\ox^*)$ is
defined by $m \mapsto \chi^m$ for $m \in N_1^\perp$.  Going the other
way, we need to study what happens to
$(\ox(D_\rho),\iota_\rho,c_{\chi^m})$ under the map
\d C_\Delta(X) \longrightarrow C_{\Delta_1}(X)\times{\rm
Hom}_\z(N_1^\perp,H^0(X,\ox^*))\ .
\d
Let's start with the second factor.  Here, note that for $m
\in N_1^\perp$, $c_{\chi^m} : \ox \simeq \ox$ is multiplication by
$\chi^m$.  Hence the induced map $N_1^\perp \to H^0(X,\ox^*)$ is
exactly the above map $\phi$.  As for the first factor, we get
$(\ox(D_\rho), \iota_\rho,c'_{m_1})$, where $c'_{m_1} =
c_{\chi^{\a(m_1)}}$ for $m_1 \in M_1$.  However, since $\pi_1 : X \to
X_1$ is a toric map taking $n_\rho \in N$ to $n_\rho \in N_1$, it
follows easily that $\pi_1^*{\cal O}_{X_1}(D_\rho) \simeq
\ox(D_\rho)$ in a way that preserves the section $\iota_\rho$ (this
follows, for example, by looking at line bundles as determined by
support functions and studying how $\pi_1$ affects support
functions).  For $m_1 \in M_1$, we have $\chi^{m_1}\circ \pi_1 =
\chi^{\a(m_1)}$, and it follows immediately that $\pi_1^*({\cal
O}_{X_1}(D_\rho), \iota_\rho,c_{\chi^{m_1}}) \sim
(\ox(D_\rho),\iota_\rho,c'_{m_1})$.  The proves commutivity, and the
theorem follows.\qed
\medskip

\noindent{\bf Remark 1.1.} When the $n_\rho$'s span $N_\r$, we get an
alternate description of the universal \xcoll\ as follows.  By [2],
$\a_\rho = [D_\rho] \in \pic$ gives a sheaf $\ox(\a_\rho)$ on $X$,
which is a line bundle since $X$ is smooth.  Furthermore, [2] gives a
canonical isomorphism $S_{\a_\rho} \simeq H^0(X,\ox(\a_\rho))$.  Since
$\deg(x_\rho) = \a_\rho$ in $\pic$, we have $x_\rho \in S_{\a_\rho}$,
so that we can write $x_\rho \in H^0(X,\ox(\a_\rho))$.  Finally, if $m
\in M$, then $\sum_\rho \langle m,n_\rho\rangle \a_\rho = 0$ in
$\pic$, which gives a canonical isomorphism
\d c_m : \otimes_\rho \ox(\a_\rho)^{\otimes \langle m,n_\rho\rangle}
\simeq \ox\big(\textstyle{\sum_\rho} \langle m,n_\rho\rangle
\a_\rho\big) = \ox\ .
\d
Then $(\ox(\a_\rho),x_\rho,c_m)$ is equivalent to the universal \xcoll\
$(\ox(D_\rho),\iota_\rho,c_{\chi^m})$.  This follows easily using the
isomorphisms $\ox(\a_\rho) \simeq \ox(D_\rho)$ constructed in [2,
\S3].
\medskip

\noindent {\bf Remark 1.2.}  Using the representability criterion
given in Proposition 4.5.4 from [3], one can prove directly that
$C_\Delta$ is representable, without knowing the toric variety $X$.
To see how this works, let $\sigma \in \Delta_{max}$, and define the
functor $C_\Delta^\sigma$ by
\d C_\Delta^\sigma(Y) = \{(L_\rho,u_\rho,c_m) \in C_\Delta(Y) :
\hbox{$u_\rho$ is an isomorphism for all $\rho \not\subset \sigma$}\}\ .
\d
Using the isomorphisms $u_\rho^{-1} : L_\rho \simeq \oy$ for $\rho
\not\subset \sigma$, one gets an equivalence $(L_\rho,u_\rho,c_m) \sim
(L'_\rho, u'_\rho, c'_m)$ where $L'_\rho = \oy$ and $u'_\rho = 1$
whenever $\rho \not\subset \sigma$.  From here, the techniques of
Examples~1.1 and 1.3 and Theorem~1.1 can be adapted to show that
$C_\Delta^\sigma$ is represented by $\A^d_k \times_k {\bb G}^{n-d}_m$,
where $d$ is the dimension of $\sigma$.  According to Proposition
4.5.4 of [3], $C_\Delta$ is then representable provided we can show
the following:
\vskip0pt
\item{(i)} The natural transformation $C_\Delta^\sigma \to C_\Delta$ is
representable by an open immersion.
\item{(ii)} The functor $C_\Delta$ is a sheaf when restricted to open
subsets of $Y$.
\item{(iii)} $C_\Delta$ is the union (as defined in part (iii) of
Proposition 4.5.4 of [3]) of the $C_\Delta^\sigma$.
\vskip0pt
\noindent The proof of (ii) is completely straightforward, and (iii)
follows easily from the nondegeneracy condition.  For (i), we need to
show that given a \xcoll\ $(L_\rho,u_\rho,c_m)$ on $Z$, the functor $Y
\mapsto \{g \in \hom(Y,Z) : g^*(L_\rho,u_\rho, c_m) \in
C_\Delta^\sigma(Y)\}$ is representable by an open subset $Z_\sigma \subset
Z$.  This is easy: $Z_\sigma$ is the biggest open subset of $Z$ where
$u_\rho$ is an isomorphism for all $\rho \not\subset \sigma$.  We
leave the details to the reader.

	By proving that $C_\Delta$ is representable, we get an
alternate construction of the smooth toric variety $X$.  This might be
useful for studying toric varieties over more general bases (for
example, over the integers or over finite fields).
\bigskip

\noindent {\bf \S2. Maps Between Toric Varieties.}  As an application
of Theorem 1.1, we will describe all maps from $\P^m_k$ to a smooth
toric variety $X$ where the $n_\rho$'s span $N_\r$.  In this case,
recall that $X$ is the geometric quotient $(\A^{\Delta(1)}_k - Z)/G$.

\proclaim Theorem 2.1. Let $X$ be a smooth toric variety such that the
$n_\rho$'s span $N_\r$, and suppose we have homogeneous polynomials
$P_\rho \in k[t_0,\dots,t_m]$, indexed by $\rho \in \Delta(1)$, such
that:
\vskip0pt
\item{(a)} If $P_\rho$ has degree $d_\rho$, then $\sum_\rho d_\rho
n_\rho = 0$ in $N$.
\item{(b)} $(P_\rho(t_0,\dots,t_m)) \notin Z$ in $\A_k^{\Delta(1)}$
whenever $(t_0,\dots,t_m) \ne 0$ in $\A^{m+1}_k$.
\vskip0pt
\noindent If we define $\tilde f(t_0,\dots,t_m) =
(P_\rho(t_0,\dots,t_m)) \in \A^{\Delta(1)}_k$, then there is a
morphism $f: \P^m_k \to X$ such that the diagram
\d \matrix{\A^{m+1}_k -\{0\} & \mapname{\tilde f} & \A_k^{\Delta(1)}-Z
\cr \downarrow && \downarrow \cr \P^m_k & \mapname{f} & X \cr}
\d
commutes, where the vertical maps are the quotient maps.  Furthermore:
\vskip0pt
\item{(i)} Two sets of polynomials $\{P_\rho\}$ and $\{P'_\rho\}$
determine the same morphism $f: \P^m_k \to X$ if and only if there is
$g \in G = {\rm Hom}_\z(\pic,{\bb G}_m)$ such that $P'_\rho =
g([D_\rho]) P_\rho$ for all $\rho$.
\item{(ii)} All morphisms $f : \P^m_k \to X$ arise in this way.
\vskip0pt

\prf Given the $P_\rho$'s satisfying (a) and (b), note that for every
$m \in M$ we have $\sum_\rho d_\rho\langle m,n_\rho\rangle = 0$, which
gives a canonical isomorphism of sheaves
\d c_m^{can} : \otimes_\rho \oP(d_\rho)^{\langle m,n_\rho\rangle}
\simeq  \oP\big(\sum_\rho d_\rho\langle m,n_\rho\rangle\big) = \oP\ .
\d
Then $(\oP(d_\rho), P_\rho,c_m^{can})$ is clearly a
$\Delta$-collection, so that we get a map $f : \P^m_k \to X$.  Using
the arguments from Theorem 1.1, one can show that if $\pi :
\A^{m+1}_k-\{0\} \to \P^m_k$ is the quotient map, then
$\pi^*(\oP(d_\rho), P_\rho,c_m^{can}) \sim ({\cal
O}_{\A^{m+1}_k-\{0\}}, P_\rho,1)$.  From here, the commutivity of the
diagram follows easily.

	Now suppose that two sets of polynomials $\{P_\rho\}$ and
$\{P'_\rho\}$ give the same map $f$.  Then, by Theorem 1.1, we know
that $(\oP(d_\rho), P_\rho,c_m^{can}) \sim
(\oP(d_\rho),P'_\rho,c_m^{can})$.  This means that there are constants
$\lambda_\rho \in k^*$ such that $P'_\rho = \lambda_\rho P_\rho$ and
$\prod_\rho \lambda_\rho^{\langle m,n_\rho\rangle} = 1$ for all $m \in
M$ because $c_m^{can}$ is preserved.  As in the proof of Theorem 1.1,
this implies that there is $g \in G$ such that $g([D_\rho]) =
\lambda_\rho$ for all $\rho$, and (i) is proved.

	Finally, to prove (ii), let $f : \P^m_k \to X$ be a morphism.
By Theorem 1.1, we know that $f$ is determined by some \xcoll\
$(L_\rho,u_\rho,c_m)$.  Since each $L_\rho \simeq \oP(d_\rho)$ for
some $d_\rho$, we get an equivalence $(L_\rho,u_\rho,c_m) \sim
(\oP(d_\rho), F_\rho,c'_m)$.  Then $(c'_m)^{-1}\circ c_m^{can} : \oP
\simeq \oP$, and thus, as in the second paragraph of the proof of
Theorem 1, we can find $\lambda_\rho \in k^*$ such that $c_m^{can} =
\prod_\rho \lambda_\rho^{\langle m,n_\rho\rangle} c'_m$ for all $m$
(this uses our assumption that the $n_\rho$'s span $N_\r$).  If we set
$P_\rho = \lambda_\rho F_\rho$, then $(\oP(d_\rho),F_\rho,c'_m) \sim
(\oP(d_\rho),P_\rho,c_m^{can})$, which shows that $f$ is determined by
the $P_\rho$'s.  It is easy to see that conditions (a) and (b) are
satisfied, and the theorem is proved.\qed
\medskip

\noindent {\bf Remark 2.1.} When $X = \P^n_k$, the $n_\rho$'s are
$e_0,\dots,e_n$ as in Example 1.2.  One can check that $\sum_{i=0}^n
d_i e_i = 0$ if and only if $d_0 = \cdots = d_n$, and then Theorem 2.1
gives the usual description of maps between $\P^m_k$ and $\P^n_k$.
\medskip

\noindent {\bf Remark 2.2.} Theorem 2.1 applies to all smooth {\it
complete\/} toric varieties since the $n_\rho$'s obviously span $N_\r$
in this case.
\medskip

\noindent {\bf Remark 2.3.} When the $n_\rho$'s don't span $N_\r$,
then, as in the proof of Theorem~1.1, we can write $X \simeq
X_1\times_k T_1$, where $T_1 \simeq {\bb G}^{n-d}_m$.  In this case, a
map $\P^m_k \to X$ is determined by maps $\P^m_k \to X_1$ and $\P^m_k
\to {\bb G}_m^{n-d}$.  The first of these maps can be described by
Theorem~2.1, and the second is obviously constant.  Thus we can
describe maps from $\P^m_k$ to an arbitrary smooth toric variety $X$.
\medskip

	We should also mention that Theorem 2.1 has been used by Guest
to study the topology of the space of rational curves on $X$ (see
[4]).

	Finally, we will discuss a more general version of Theorem 2.1,
where $\P^m_k$ is replaced by an arbitrary complete toric variety
$Y$.  If $Y$ is determined by the fan $\Delta_Y$, then by [2], $Y$ has
a homogeneous coordinate ring $S^Y = k[y_\tau]$, where $\tau \in
\Delta_Y(1)$.  The ring $S^Y$ is graded by the Chow group
$A_{n-1}(Y)$, and we denote the graded pieces by $S^Y_\a$ for $\a \in
A_{n-1}(Y)$.  Note also that ${\rm Pic}(Y) \subset A_{n-1}(Y)$.  By
[2], we can also express $Y$ as a categorical quotient of
$\A^{\Delta_Y(1)}_k - Z_Y$.  Then we get the following result:

\proclaim Theorem 2.2. Let $X$ be a smooth toric variety such that the
$n_\rho$'s span $N_\r$, and let $Y$ be a complete toric variety with
coordinate ring $S^Y$.  Suppose we have homogeneous polynomials
$P_\rho \in S^Y$, indexed by $\rho \in \Delta(1)$, such
that:
\vskip0pt
\item{(a)} If $P_\rho \in S^Y_{\beta_\rho}$, then $\beta_\rho \in {\rm
Pic}(Y)$ and $\sum_\rho \beta_\rho\otimes n_\rho = 0$ in ${\rm
Pic}(Y)\otimes N$.
\item{(b)} $(P_\rho(t_\tau)) \notin Z$ in $\A_k^{\Delta(1)}$
whenever $(t_\tau) \notin Z_Y$ in $\A^{\Delta_Y(1)}_k$.
\vskip0pt
\noindent If we define $\tilde f(t_\tau) = (P_\rho(t_\tau)) \in
\A^{\Delta(1)}_k$, then there is a
morphism $f: Y \to X$ such that the diagram
\d \matrix{\A^{\Delta_Y(1)}_k - Z_Y & \mapname{\tilde f} &
\A_k^{\Delta(1)}-Z \cr \downarrow && \downarrow \cr Y & \mapname{f} &
X \cr}
\d
commutes, where the vertical maps are the quotient maps.  Furthermore:
\vskip0pt
\item{(i)} Two sets of polynomials $\{P_\rho\}$ and $\{P'_\rho\}$
determine the same morphism $f: Y \to X$ if and only if there is
$g \in G = {\rm Hom}_\z(\pic,{\bb G}_m)$ such that $P'_\rho =
g([D_\rho]) P_\rho$ for all $\rho$.
\item{(ii)} All morphisms $f : Y \to X$ arise in this way.
\vskip0pt

\prf We will only sketch the proof, leaving the details to the reader.
The key idea is that by [2], $\a \in {\rm Pic}(Y)$ gives a line bundle
$\oy(\a)$ such that we have canonical isomorphisms $\oy(\a)\otimes
\oy(\beta) \to \oy(\a+\beta)$ for $\a,\beta \in {\rm Pic}(Y)$.
Furthermore, from [2] there is a natural isomorphism $H^0(Y,\oy(\a))
\simeq S^Y_\a$.  Then it is easy to see that the $P_\rho$'s give a
\xcoll\ $(\oy(\a_\rho),P_\rho,c^{can}_m)$, and from here the rest
of the proof is identical to what we did in Theorem 2.1.\qed
\bigskip

\noindent {\bf \S3. Concluding Remarks. } Another description of the
functor represented by a toric variety $X$ is due to Ash, Mumford,
Rapoport and Tai (see [1, Chapter I, \S2]).  They consider pairs
$({\cal S},\pi)$ such that:
\vskip0pt
\item{(i)} ${\cal S}$ is a sheaf of sub-semigroups of the constant sheaf
$M_Y$ on $Y$ determined by $M$.
\item{(ii)} $\pi: {\cal S} \to \oy$ is a semigroup homomorphism ($\oy$
is a semigroup under multiplication).
\vskip0pt
\noindent Furthermore, they assume that $({\cal S},\pi)$ has the
following properties:
\vskip0pt
\item{(iii)} For $s \in {\cal S}$, $\pi(s)$ is invertible if and only
if $s$ is.
\item{(iv)} For each $y \in Y$, there is some $\sigma \in \Delta$ such
that ${\cal S}_y = \sigma^\vee\cap M$.
\vskip0pt
\noindent The main result of [1, Chapter I, \S2] is that for all $Y$,
there is a natural bijection
\d \hom(Y,X) \simeq \{\hbox{all pairs $({\cal S},\pi)$ on $Y$
satisfying (i)--(iv) above}\}\ .
\d
This description of the functor represented by $X$ is clearly related
to the usual way of constructing $X$ by patching together the affine
schemes $X_\sigma = {\rm Spec}(k[\sigma^\vee\cap M])$.

	In contrast, our description of $\hom(Y,X)$ is more closely
tied to the geometric quotient $X \simeq (\A_k^{\Delta(1)} - Z)/G$.
An advantage of our approach is how it generalizes the usual
description of maps between projective spaces (see Theorem 2.1).  The
Ash-Mumford-Rapoport-Tai approach, on the other hand, has the virtue
that it applies to all toric varieties, not just smooth ones.  (The
problem with our description in the nonsmooth case is that
$\ox(D_\rho)$ need not be a line bundle.)  It would be interesting to
see the analog of Theorems 1.1 and 2.1 for the case of simplicial
toric varieties.

	I am grateful to Martin Guest and Stein Arild Str\o mme for
bringing this problem to my attention.  The research for this paper
was supported by NSF grant DMS-9301161.
\bigskip

\noindent {\bf References}
\medskip

\item{1.} A. Ash, D. Mumford, M. Rapoport and Y.-S. Tai, {\sl Smooth
Compactifications of Locally Symmetric Varieties\/}, Math Sci Press,
Brookline, MA, 1975.
\medskip

\item{2.} D. Cox, {\it The homogeneous coordinate ring of a toric
variety\/}, to appear.
\medskip

\item{3.} J. Dieudonn\'e and A. Grothendieck, {\sl El\'ements de
G\'eom\'etrie Alg\'ebrique I\/}, Springer-Verlag, Berlin Heidelberg
New York, 1971.
\medskip

\item{4.} M. Guest, {\it The topology of the space of rational curves
on a toric variety\/}, preprint.
\medskip

\end